\documentclass[12pt]{article}

\usepackage{amssymb}

\date{}
\title{\vspace{-1cm}A note on embedding hypertrees}
\author{
Po-Shen Loh \thanks{Research supported in part by a Fannie and John Hertz Foundation Fellowship, an NSF Graduate
Research Fellowship, and a Princeton Centennial Fellowship.}\\
\small Department of Mathematics\\[-0.8ex]
\small Princeton University\\
\small \texttt{ploh@math.princeton.edu}
}

\newtheorem{theorem}{Theorem}

\begin{document}
\maketitle

\begin{abstract}
  A classical result from graph theory is that every graph with
  chromatic number $\chi > t$ contains a subgraph with all degrees at
  least $t$, and therefore contains a copy of every $t$-edge tree.
  Bohman, Frieze, and Mubayi recently posed this problem for
  $r$-uniform hypergraphs.  An \emph{$r$-tree}\/ is a connected
  $r$-uniform hypergraph with no pair of edges intersecting in more
  than one vertex, and no sequence of distinct vertices and edges
  $(v_1, e_1, \ldots, v_k, e_k)$ with all $e_i \ni \{v_i, v_{i+1}\}$,
  where we take $v_{k+1}$ to be $v_1$.  Bohman, Frieze, and Mubayi
  proved that $\chi > 2rt$ is sufficient to embed every $r$-tree with
  $t$ edges, and asked whether the dependence on $r$ was necessary.
  In this note, we completely solve their problem, proving the tight
  result that $\chi > t$ is sufficient to embed any $r$-tree with $t$
  edges.
\end{abstract}

\section{Introduction}

An \emph{$r$-graph}\/ is a hypergraph where all edges have size $r$,
and a proper coloring is an assignment of a color to each vertex such
that no edge is monochromatic.  The \emph{chromatic number $\chi$} is
the minimum $k$ for which there is a proper coloring with $k$ colors.
A natural question is to investigate what properties can be forced by
sufficiently large chromatic number.  In the case of graphs, much is
known, from trivialities such as $\chi > t$ implying the existence of
a subgraph with all degrees at least $t$, to deeper results such as
$\chi > 4$ implying non-planarity.  Far less is known for hypergraphs,
but a folklore observation (see, e.g., \cite{EL}) is that whenever
$\chi > 2$, there is a pair of edges that intersect in a single
vertex.  This structure corresponds to a 2-edge \emph{hypertree},
which in general is a connected hypergraph with no pair of edges
intersecting in more than one vertex, and no sequence of distinct
vertices and edges $(v_1, e_1, \ldots, v_k, e_k)$ with all $e_i \ni
\{v_i, v_{i+1}\}$, where we take $v_{k+1}$ to be $v_1$.

For graphs, $\chi > t$ implies that there is a subgraph with all
degrees at least $t$, in which we can embed any $t$-edge tree.
Bohman, Frieze, and Mubayi recently posed the problem of generalizing
this result to $r$-graphs.  As they noted, this is not entirely
trivial because there are hypergraphs with arbitrarily large minimum
degree, but no copy of the path with 3 edges.  Indeed, consider the
3-graph with vertex set $\{v_1, \ldots, v_n\}$ and edges consisting of
all triples containing $v_1$.

Observe that an $r$-uniform hypertree (henceforth referred to as an
\emph{$r$-tree}) with $t$ edges always has exactly $1 + (r-1)t$
vertices.  So, the complete $r$-graph on $(r-1)t$ vertices does not
contain any $r$-tree with $t$ edges, while its chromatic number is
exactly $t$.  On the other hand, Bohman, Frieze, and Mubayi proved in
\cite{BFM} that every $r$-graph with $\chi > 2rt$ contains a copy of
every $r$-tree with $t$ edges.  They believed that their bound was far
from the truth, and remarked at the end of their paper that it would
be interesting to determine whether it should depend on $r$ in an
essential way.  In this note, we completely solve their problem,
proving the following tight result.

\begin{theorem}
  \label{thm:main}
  Every $r$-uniform hypergraph with chromatic number greater than $t$
  contains a copy of every $r$-uniform hypertree with $t$ edges.
\end{theorem}

\section{Proof}

It suffices to show that for any $r$-tree $T$ with $t$ edges, every
$T$-free $r$-graph $H$ can be properly colored with the integers $\{1,
\ldots, t\}$.  Although the proof is short, the following special case
helps to illuminate the argument.  Suppose the $r$-tree $T$ is a path
with $t$ edges, and there is a proper $t$-coloring of $H-e_1^*$, the
hypergraph on the same vertex set but with an arbitrary edge $e_1^*$
removed.  The edge $e_1^*$ is monochromatic, say in color 1, or else
we are done.  Let $v_1^*$ be an arbitrary vertex of $e_1^*$.  Either
we can recolor $v_1^*$ in color 2 without making any edge
monochromatic in color 2 (and hence are done because $e_1^*$ is no
longer monochromatic), or else some edge $e_2^* \ni v_1^*$ has all
vertices except $v_1^*$ colored 2.  Note that since all vertices in
$e_2^*$ are colored 2 except for $v_1^*$, and all vertices in $e_1^*$
are colored 1, the two edges intersect only at $v_1^*$, thus forming a
copy of the 2-edge path.

Suppose for a moment that $e_2^*$ is the unique edge containing
$v_1^*$ which has all vertices except $v_1^*$ colored 2.  Repeating
the argument, we select $v_2^* \in e_2^*$, and either find an edge
$e_3^* \ni v_2^*$ with all other vertices colored 3 (thus forming a
3-edge path together with $e_2^*$ and $e_1^*$), or obtain a proper
coloring of $H$ by recoloring $v_2^*$ with color 3 and $v_1^*$ with
color 2.  Unfortunately, when $e_2^*$ is not unique, the recoloring of
$v_1^*$ with color 2 may make another edge monochromatic, so a more
careful argument is needed in general.  Nevertheless, for illustration
only, let us make the simplifying uniqueness assumption, and continue
in this way to find successively longer paths $e_1^*, e_2^*, \ldots,
e_s^*$.  Yet $H$ has no $t$-edge path, so this must stop before we
need to use $t+1$ colors.  Then, we will be able to properly $t$-color
$H$ by recoloring each vertex $v_i^*$ with color $i+1$.

\vspace{3mm}

\noindent \textbf{Proof of Theorem 1.}\, Let $T$ be an $r$-tree with
$t$ edges.  We will show that every $T$-free $r$-graph $H$ can be
properly colored with the integers $\{1, \ldots, t\}$.  Preprocess $T$
by labeling its edges and coloring its vertices as follows.  Let $e_1$
be an arbitrary edge of $T$, and label the other edges with $e_2,
\ldots, e_t$ such that for each $i \geq 2$, all edges $e_j$ along the
(unique) path linking $e_i$ and $e_1$ are indexed with $j < i$.  This
can be done by exploring $T$ via breadth-first-search, for instance.
Then, color each vertex $v \in T$ with the integer equal to the
minimal index $i$ for which $e_i \ni v$.

We now induct on the number of edges of $H$.  Let $e_1^*$ be an edge
of $H$, and suppose that there is a proper $t$-coloring of $H -
e_1^*$, the hypergraph on the same vertex set, but without the edge
$e_1^*$.  If this is already a proper coloring of $H$, then we are
done.  Otherwise, without loss of generality all vertices of $e_1^*$
received the color 1.  The following recoloring algorithm formalizes
the above heuristic.
\begin{enumerate}
\item Let $H' \subset H$ be a maximal colored-copy of a subtree of $T$
  containing $e_1$, and let $T' \subset T$ be that subtree.  This
  means there is a color-preserving injective graph homomorphism $\phi
  : T' \rightarrow H$ with maximal $T' \ni e_1$, which exists because
  $e_1^*$ itself is a colored-copy of $e_1$.

\item Since $H$ is $T$-free, there is an edge $e_s$ in $T$ but not
  $T'$, which is incident to some vertex $v \in T'$.  Change the color
  of $\phi(v) \in H$ to $s$.  Terminate if $\phi(v) \in e_1^*$;
  otherwise, return to step 1.
\end{enumerate}
The maximality of $H'$ ensures that the recoloring step never creates
any new monochromatic edges.  Indeed, suppose for contradiction that
$H$ has an edge $e' \ni \phi(v)$ with all vertices except $\phi(v)$
colored $s$.  Our preprocessing of $T$ ensures that no vertex in the
colored-copy $H'$ of $T'$ has color $s$, so $e'$ intersects $H'$ only
at $\phi(v)$.  Thus $H' + e'$ would be a colored copy of $T' + e_s$,
contradicting maximality.

Also, the algorithm terminates because the recoloring step always
increases the (integer) color of $\phi(v)$, but no color ever exceeds
$t$.  To see this, observe that since we had a colored-copy, the color
of $\phi(v)$ originally equalled the color of $v \in T$, which we
defined to be the minimal index $i$ such that $e_i \ni v$.  By our
preprocessing of $T$, $e_s \not \in T'$ implies that some
lower-indexed edge also contains $v$.  Hence $\phi(v)$ indeed had
color less than $s$.  Therefore, we eventually obtain a proper
coloring of $H$.  \hfill $\Box$

\section{Concluding remarks}

\begin{itemize}
\item The standard proof of the graph case of Theorem \ref{thm:main}
  uses the fact that every $t$-edge tree can be embedded in any graph
  with minimum degree at least $t$.  This is not true for hypergraphs,
  so our proof uses a completely different argument that does not rely
  on degrees at all.  Consequently, our proof also gives a new
  perspective on the graph case.

\item Results for graphs that used $\chi > t$ to embed $t$-edge trees
  can now be extended to uniform hypergraphs.  Consider, for example,
  the following classical result of Chv\'atal, referred to as ``one of
  the most elegant results of Graph Ramsey Theory'' by Graham,
  Rothschild, and Spencer in their book \cite{GRS}.  The Graph Ramsey
  number $R(H_1, H_2)$ is the smallest $n$ such that every red-blue
  edge-coloring of $K_n$ contains either a red copy of $H_1$ or a blue
  copy of $H_2$.  When $H_1$ is a complete graph $K_k$ and $H_2$ is
  any $t$-edge tree, Chv\'atal determined that $R(H_1, H_2)$ is
  precisely $(k-1)t + 1$.

  Using Theorem \ref{thm:main}, we can lift one of the standard proofs
  of this result to $r$-graphs.  Indeed, suppose we have a red-blue
  edge-coloring of the complete $r$-graph on $(k-1)t + 1$ vertices,
  and let $H$ be the hypergraph on the same vertex set formed by
  taking only the blue edges.  If $\chi(H) \leq t$, then $H$ has an
  independent set of size at least $\big\lceil \frac{(k-1)t + 1}{t}
  \big\rceil = k$, which corresponds to a red complete $r$-graph on
  that many vertices.  Otherwise, if $\chi(H) > t$, then Theorem
  \ref{thm:main} implies that any $t$-edge tree can be found in the
  blue graph $H$.
  
  On the other hand, the $r$-graph obtained by taking the disjoint
  union of $\big\lfloor \frac{k-1}{r-1} \big\rfloor$ copies of the
  complete $r$-graph on $(r-1)t$ vertices does not contain any
  $r$-tree with $t$ edges, while its independence number is at most
  $k-1$.  So, if we color all of its edges blue, and add in all
  missing edges with color red, then we obtain an edge-coloring of the
  complete $r$-graph on $\big\lfloor \frac{k-1}{r-1} \big\rfloor \cdot
  (r-1)t$ vertices with no red $K_k^{(r)}$ and no blue $r$-tree with
  $t$ edges.  Therefore, the $r$-graph result is tight when $r-1$
  divides $k-1$, and asymptotically tight for $k \gg r$.
\end{itemize}

\noindent {\bf Acknowledgment.}\, The author thanks his Ph.D. advisor,
Benny Sudakov, for introducing him to this problem, and for remarks
that helped to improve the exposition of this note.  Also, he thanks
Asaf Shapira for pointing out the application of the main theorem to
Chv\'atal's result, and the referee for carefully reading this
article.

\end{document}